%% file: JOC-template.tex
\documentclass[ijoc,nonblindrev]{informs3} 

\OneAndAHalfSpacedXII 

\usepackage{natbib}
 \bibpunct[, ]{(}{)}{,}{a}{}{,}%
 \def\bibfont{\small}%
\TheoremsNumberedThrough     

\EquationsNumberedThrough    

\usepackage{rotating}
\usepackage[theme=grayscale]{jlcode}
\usepackage[utf8]{inputenc}

\usepackage{hyperref}
\usepackage[anythingbreaks]{breakurl}

\usepackage{array,multirow,graphicx}
\usepackage{float}
\usepackage{tabularx,dcolumn,booktabs}

\DeclareMathOperator{\st}{s.t.}

\begin{document}


\RUNAUTHOR{Dias Garcia, Bodin and Street}

\RUNTITLE{\texttt{BilevelJuMP.jl}: Modeling and Solving Bilevel Optimization in Julia}

\TITLE{\texttt{BilevelJuMP.jl}:  Modeling and Solving Bilevel Optimization in Julia}

\ARTICLEAUTHORS{%
\AUTHOR{Joaquim Dias Garcia}
\AFF{LAMPS at PUC-Rio \& PSR, \EMAIL{joaquim@psr-inc.com}}
\AUTHOR{Guilherme Bodin}
\AFF{LAMPS at PUC-Rio \& PSR, \EMAIL{guilhermebodin@psr-inc.com}}
\AUTHOR{Alexandre Street}
\AFF{LAMPS at PUC-Rio, \EMAIL{street@ele.puc-rio.br}}
} 

\ABSTRACT{%
In this paper we present BilevelJuMP, a new Julia package to support bilevel optimization within the JuMP framework. The package is a Julia library that enables the user to describe both upper and lower-level optimization problems using the JuMP algebraic syntax. Due to the generality and flexibility our library inherits from JuMP's syntax, our package allows users to model bilevel optimization problems with conic constraints in the lower level and all JuMP supported constraints in the upper level (Conic, Quadratic, Non-Linear, Integer, etc.). Moreover, the user-defined problem can be subsequently solved by various techniques relying on mathematical program with equilibrium constraints (MPEC) reformulations. Manipulations on the original problem data are possible due to MathOptInterface.jl's structures and Dualization.jl features. Hence, the proposed package allows quickly model, deploy, and thereby experiment bilevel models based on off-the-shelf mixed integer linear programming and nonlinear solvers.
}%



\maketitle

%


\section{Introduction}\label{intro}

Bilevel optimization has been a widely used modeling tool in mathematical programming, operations research and economics since its first introduction by \cite{von1952theory} in game theory. The broad range of applications include hyperparameter optimization in machine learning \citep{franceschi2018bilevel}, toll setting in transportation networks \citep{labbe1998bilevel}, multiple problems in energy and power systems \citep{pozo2017basic}, defense applications \citep{bracken1974defense}, facility location \citep{kuccukaydin2011competitive} only to list a few areas.

Complete introductions to bilevel optimization can be found on books covering theoretical background and analysis, taxonomy, solutions algorithms for special classes, and selected applications \citep{dempe2002foundations,bard2013practical,dempe2015bilevel}.
In addition to books, many reviews on the subject have been published in the last decades \citep{vicente1994bilevel,colson2007overview,kalashnikov2015bilevel,beck2021gentle}.
A very long list of publications related to bilevel optimization can be found in \cite{dempe2018bilevel}.

It is well known that general bilevel optimization problems fall in the NP-hard class \citep{jeroslow1985polynomial}. Hence, there is no hope in finding efficient algorithms for generic problems. On the other hand, modeling bilevel problems is a theoretically more straightforward task, albeit, in practice, the modeling step can be the difference between finding a tractable model for which there is a reasonable solution approach in realistic cases or not. Because bilevel models are very complex and constitute a broad class of mathematical programming problems, many modeling languages lack the proper functionality to handle these problems.

In the following subsections, we present the main literature regarding available techniques and software to place and justify the contribution of our work properly. Notwithstanding, it is out of the scope of this section to provide a comprehensive review on the subject, for which we refer to previously reported books and reviews.

\subsection{Solving bilevel optimization}

Many strategies have been proposed to solve Bilevel Optimization problems. Some of the most widely known techniques are based on classical algorithms such as the simplex method, the branch and bound method, and the interior-point methods.
Due to the inherent combinatorial nature of many Bilevel optimization problems, some of the early developed techniques are intrinsically combinatorial. Among a large set of enumeration algorithms, which can be seen as modifications of the simplex method for linear programming, we highlight the $K^{th}$-best algorithm \citep{bialas1982two}.

The fundamental technique that will be explored in this work is converting the bilevel problem into a single-level problem by adding the lower-level KKT conditions to the upper-level problem. The resulting optimization problem is known as Mathematical Programming with Equilibrium Constraints, MPEC \citep{kim2020mpec}. This group of techniques has also been labeled enumeration-based because the main difficulty is handling the complementarity constraints. Thus, a classic solution is a specially tailored Branch and Bound \citep{bard1988convex,hansen1992new}. Instead of writing a branch and bound method from scratch, one could reformulate the single-level problem into an amenable form for an off-the-shelf algorithm like mixed integer programming (MIP) solvers. This method was first presented relying upon big-M formulations in \cite{fortuny1981representation}. More recently, a special ordered set of type 1, SOS1, formulations were developed in \cite{siddiqui2013sos1}.

Interestingly, non-linear programming problem (NLP), with additional regularization terms, can also be used to solve MPECs \citep{fletcher2004solving,scholtes2001convergence,ralph2004some}. In this case, it would not make as much sense to call the MPEC reformulations an enumeration-based method. On the other hand, the latter might be called local search methods as this strategy leads to local solutions, in contrast with the global ones provided by MIP-based methods.
Some solvers were specially tailored to tackle MPEC's: KNITRO \citep{nocedal2006knitro}, filterMPEC \citep{fletcher2021filterMPEC, fletcher2004solving}, NLPEC \citep{gams2021nlpec}. A combination of NLP and MIP based methods was proposed by \citep{pineda2018efficiently}.

Other strategies to solve bilevel problems include: bundle type algorithms \citep{falk1995bilevel}, semi-definite relaxations \citep{fampa2013semidefinite}, penalty function based methods \citep{aiyoshi1984solution,white1993penalty,kleinert2020computing}, Benders decomposition \citep{byeon2019benders}.
Cutting-planes approaches \citep{wu1998cutting, fischetti2017new, tahernejad2020branch} have received attention recently because they can handle lower-level problems with integer variables --- solvers were developed by \cite{ted_ralphs_2019_3560359, fischetti2020solver}.
Descent methods were proposed by \cite{kolstad1990derivative} and \cite{vicente1994descent} to obtain quick local solutions. Heuristic methods were also developed to obtain practical solutions, for instance, the bi-objective-based method of \cite{bard1983efficient}.
Finally, we refer to \cite{sinha2017review} for evolutionary approaches.

\subsection{Modeling bilevel optimization}

Algebraic modeling languages (AML) play a central role in unlocking the huge potential of optimization models through a mathematical friendly environment that integrates solvers and models in a very practical manner where problems of all disciplines can be efficiently modeled and solved. Bilevel Optimization interfaces, which mostly automate reformulations and pass to specific solver types, have also been proposed.

GAMS has an interface described in its Extended Mathematical Programming, \cite{GAMS2021bilevel}. Variables and constraints are created as usual, and then they are annotated to specify the level they belong to in an external file. The annotated problem is reformulated by GAMS EMP routines using the KKT reformulations \citep{ferris2009extended}. Finally, the problem is optimized by the available MPEC solvers, namely, the above cited KNITRO and NLPEC. The follower subproblems can be linearly constrained with quadratic objectives (QP) or Variational Inequalities (VI). The upper-level problem is constrained by the selected solver capability.

The YALMIP MATLAB package for optimization modeling \citep{lofberg2004yalmip} also provides a Bilevel Optimization interface \cite{YALMIP2021bilevel}. The variables and constraints are defined by the standard methods, and then they are passed as lists to a ``solvebilevel" function. The lower level problem can be a QP, and the upper level can be anything supported by the YALMIP interface. The available solution methods are based on the MPEC reformulation where the handling of complementarity constraints is forwarded to external MIP, NLP or MPEC solvers or to an internal branch and bound that allows different solvers for upper and lower-level problems.

It is also possible to model Bilevel optimization in the Pyomo Adversarial Optimization, PAO  \citep{hart2019pao}, a Python package extending Pyomo \citep{hart2017pyomo}. A (lower) submodel is created as an object attached to the main (upper) model, then variables and constraints are created directly in their owner models, the former can be shared in objectives and constraints. As in the previous two AMLs, the problem is automatically reformulated, and in Pyomo's case can be passed either to an NLP solver \citep{ralph2004some} or a MIP solver \citep{fortuny1981representation}. Additionally, PAO has an interface to the MibS solver \citep{ted_ralphs_2019_3560359} and implements a Column Constraint Generation, CCG \citep{yue2019projection}.

The Julia \citep{bezanson2017julia} package BilevelOptimization.jl \citep{Besancon2019BilevelJulia} provides a very simple interface for Bilevel modeling in JuMP \citep{LubinDunningIJOC,DunningHuchetteLubin2017}. Still, the interface remains very basic. For instance, although the upper level can represent arbitrary JuMP problems (as long as the selected solver supports them), the lower level is constrained to QP. Te package supports two MIP-based methods: \cite{fortuny1981representation} and the SOS1 based reformulation. The critical issue is that the lower level is not represented by a JuMP-based syntax, preventing Julia users to fully enjoy the modeling power provided by JuMP. Instead, lower-level problems must be described by matrices, which can be easily manipulated to write KKT reformulations. This is one of the salient features and the primary goal of our propose package, BilevelJuMP, namely, to allow representing the lower level problem within the JuMP syntax in a single and integrated new bilevel JuMP model.

It is important to highlight that, just like \cite{AIMMS2021comp} and \cite{AMPL2021comp}, JuMP also has native support to complementarity constraints, thereby being capable of handling MPEC models. However, none of them can model Bilevel Optimization models directly.

\subsection{Objective and Contributions}

The main objective of this work is to provide a complete open-source interface for Bilevel Optimization fully integrated in the JuMP modeling language named BilevelJuMP.jl and available online at:
\begin{center}
\url{https://github.com/joaquimg/BilevelJuMP.jl}
\end{center}
It is also readily available at the Julia Package manager, at the time of this paper, in version 0.5.1. Julia users can run \texttt{add BilevelJuMP} and have full access to all features of the library.

Regarding similar works, \cite{Besancon2019BilevelJulia} provided a great motivation and a nice first step to tackle bilevel models in Julia. However, it is incomplete as a modeling framework due to the strong limitations associated with the modeling of second-level problems as explained before. The consideration of a generic JuMP-based second-level model within a JuMP integrated interface significantly increases functionality parity with other AML's. Notwithstanding, and more importantly, it paves the way for new developments and computational applications based on bilevel optimization. Because the proposed package provides a simple an integrated interface, fully embedded into the JuMP language, new and expert users can easily prototype and test bilevel models enjoying all open-source and commercial linear, nonlinear, and MILP solvers integrated in JuMP, depending on the model characteristics. Hence, the newly proposed BilevelJuMP.jl can be used in teaching environments to introduce practical aspects of bilevel modeling as well as in practical applications inheriting many of the functionalities and advances directed to JuMP, one of the main packages of Julia.

Similar to other interfaces, BilevelJuMP.jl is also capable of reformulating bilevel problems and export the model to existing external solvers. In particular, many reformulations were implemented to allow practitioners to test each method in their particular problems. Some experiments are presented in this work to provide a first glance on the differences between existing methods.

Researchers of the Bilevel Optimization community are similarly benefited. The functionality exposed by BilevelJuMP.jl can be used as a benchmark in terms of performance and solution quality for new algorithms. These new algorithms might even be implemented in Julia with the data structures already defined in the package. Therefore, we also provide a software contribution with pieces of code that can be directly used in future implementations. 

BilevelJuMP.jl was designed to be extensible, and the various implemented methods vouch for it. As we shall discuss, the key ingredient to developing all the solution methods is to rely on MathOptInterface.jl's API.

Advanced functionality is also part of the contributions. BilevelJuMP.jl can represent a Conic Program, CP, in the lower-level problem and can deal with upper-level constraints, including dual variables of lower-level problems (see Appendix \ref{ap-dual}).
Finally, the composability inherent in many Julia packages allows performing reformulations that enable the user to solve even more complex problem classes.

\section{Conic Bilevel standard form}\label{form}

We will only consider optimistic bilevel problems \citep{dempe2002foundations}, in short, the solution of the lower level will be the one that optimizes the upper level in case of degeneracy.
We start by describing the main notation that will be used in the remainder of the work. $z$ and $x$ are vectors of decision variables, respectively, from the upper and lower-level problems. While $x$ is $n^U$ dimensional, $z$ has $n^L$ entries. $[x,z]$ is a $(n^L+n^U)$--vector with the elements of $x$ and $z$ stacked. $Q^j$, $a_i^j$, $d_i^j$, $b_i^j$, $A_i^j$, $D_i^j$, for $j\in \{U,L\}$ are matrices (upper case) and vectors (lower case) of constants. $\mathcal{C}_i^j$, for $j\in \{U,L\}$, are sets, of arbitrary finite dimension, which most commonly will be convex cones. $m^U$ and $m^L$ are the number of vector constraints in the upper and lower problems.
As in traditional bilevel programming, $z$ is decided in the upper level and passed to the lower level as a parameter and $x$ might be seen as an upper-level variable constrained to be an optimal solution of the lower level. Hence, the optimistic bilevel problem follows: 
\[
\begin{aligned}
    &\min_{x \in {\mathbb{R}^n}^L, z \in {\mathbb{R}^n}^U} && \frac{1}{2}{[x,z]}^{\top} Q^U {[x,z]} + {a_0^U}^{\top} x + {d_0^U}^{\top} z + b_0^U \\
    &\st && A_i^U x + D_i^U z +b_i^U  \in \mathcal{C}_i^U,  i = 1 \ldots m^U \\
        &&& x(z) \in
     \begin{aligned}[t]
        &\argmin_{x \in {\mathbb{R}^n}^L} && \frac{1}{2}{[x,z]}^{\top} Q^L {[x,z]} + {a_0^L}^{\top} x+ {d_0^L}^{\top} z + b_0^L\\
            &\st && A_i^L x + D_i^L z + b_i^L  \in \mathcal{C}_i^L,  i = 1 \ldots m^L \\
     \end{aligned}\\
\end{aligned}
\]

As detailed by \cite{legat2020mathoptinterface}, describing constraints as function-set pairs is very flexible. For simplicity, we limited ourselves to affine functions contained in sets in the constraints of the above model. If all sets are all convex cones, we have a standard conic form for bilevel programs.

Keeping the lower level problem as a convex conic program is especially useful for writing KKT conditions when converting the problem into MPEC form. Although lower-level integer variables could be tackled by specialized solvers \citep{tahernejad2020branch}, the same goes for non-linear constraints like the ones in \citep{sinha2017review}. The upper-level problem can be more complex, including non-linear constraints and integer variables, because they are not affected in MPEC reformulations.

\subsection{KKT Reformulation of bilevel programs}\label{kkt}

Given a conic bilevel program in the standard form that we described in the previous section, we can formulate an equivalent MPEC applying the KKT reformulation to convert the lower level optimization problem into a set of linear and non-linear equations.

Let us focus on the following parametric convex quadratic conic problem that is equivalent to the lower level problem:
\begin{align}
& \min_{x \in \mathbb{R}^n} &  + \frac{1}{2} x^{\top} Q_1 x + x^{\top} Q_2 z + \frac{1}{2} z^{\top} Q_3 z + {a_0}^{\top} x + b_0 + d_0^{\top} z \notag
\\
& \;\;\text{s.t.} & A_i x + b_i + D_i z & \in \mathcal{C}_i & i = 1 \ldots m \label{ep:primal}
\end{align}
Note that the $L$ superscripts were dropped for simplicity, and we split the $Q$ matrix in $Q_1$, $Q_2$, and $Q_3$. Because $z$ are parameters, only $Q_1$ is required to be a positive semi-definite matrix. In the following, we will denote the dual cone of $\mathcal{C}_i$ as $\mathcal{C}_i^*$.

Following \cite{boyd2004convex}, we can write the KKT conditions as:
\begin{itemize}
\item Primal Feasibility:
\begin{align}
A_i x + b_i + D_i z \in \mathcal{C}_i , \ \ i = 1 \ldots m \label{eq:kkt_primal} 
\end{align}
\item Dual Feasibility:
\begin{align}
y_i \in \mathcal{C}_i^*, \ \ i = 1 \ldots m \label{eq:kkt_dual}
\end{align}
\item Stationarity:
\begin{align}
Q_1 x + Q_2 z + a_0 - \sum_{i=1}^m A_i^{\top} y_i  = 0 \label{eq:kkt_stat}
\end{align}
\item Complementary slackness:
\begin{align}
y_i^{\top} (A_i x + b_i + D_i z) = 0, \ \ i = 1 \ldots m \label{eq:kkt_comp}
\end{align}
\end{itemize}
Putting all together we arrive at the MPEC form of the bilevel conic program:
\[
\begin{aligned}
    &\min_{x \in {\mathbb{R}^n}^L, z \in {\mathbb{R}^n}^U} && \frac{1}{2}{[x,z]}^\top Q^U {[x,z]} + {a_0^U}^\top x + {d_0^U}^\top z + b_0^U \\
    &\st && A_i^U x + D_i^U z +b_i^U  \in \mathcal{C}_i^U, \ \ i = 1 \ldots m^U \\
&&& A_i^L x + b_i^L + D_i^L z \in \mathcal{C}_i^L, \ \ i = 1 \ldots m^L \\
&&& y_i \in {\mathcal{C}_i^L}^*, \ \ i = 1 \ldots m^L \\
&&& Q_1^L x + Q_2^L z + a_0^L - \sum_{i=1}^m {A_i^L}^{\top} y_i  = 0\\
&&& y_i^\top (A_i^L x + b_i^L + D_i^L z) = 0, \ \ i = 1 \ldots m^L
\end{aligned}
\]
Such form is particularly useful to solve bilevel optimization problems as described in the previous sections. The main challenge being \eqref{eq:kkt_comp} constraints which are highly non-linear and non-convex. Dealing with the latter requires special solvers, tailor-made algorithms, or extra reformulation steps to reach the standard form of some NLP or MIP solver.

\section{BilevelJuMP.jl}

BilevelJuMP.jl is an extension of the JuMP modeling language \citep{DunningHuchetteLubin2017,LubinDunningIJOC} for optimization problems in the Julia language \citep{bezanson2017julia}. Other packages successfully extending JuMP and MOI include SDDP.jl \citep{dowson2021sddp}, SumOfSquares.jl \citep{benoit_legat_2021_4708982}, InfiniteOpt.jl \citep{pulsipher2021infiniteopt}.
BilevelJuMP.jl has two main functionalities: modeling and solving bilevel optimization problems.

This open-source package heavily relies on MathOptInterface.jl, also referred to as MOI, \citep{legat2020mathoptinterface}, another Julia package that was written to be the new backend of JuMP which led to a complete rewrite of the latter.
MathOptInterface.jl is an intermediary layer between JuMP's user-friendly AML interface and the diverse and typically matrix-oriented format of solvers.
In BilevelJuMP.jl, MOI is used to store problem data from an extended JuMP interface and reformulate bilevel optimization problems into the MPEC form and then into a solver-compatible formulation of MPEC.

\subsection{A Modeling Interface for Bilevel Optimization}

The basic modeling interface of BilvelJuMP.jl relies on JuMP's extensible methods and macros to write and combine two optimization problems. Not surprisingly, other methods had to be created to accommodate the needs of bilevel optimization interfaces.

The main data structure in this software is the \texttt{BilevelModel}, which is a subtype of JuMP's \texttt{AbstractModel}. \texttt{BilevelModel} holds two other JuMP Models to represent the upper and lower optimization problems. Also, additional information is held to link the two problems and store additional JuMP data and attributes used in reformulations.

Just like a regular JuMP \texttt{Model}, the \texttt{BilevelModel} will need a solver constructor to solve an optimization problem. On the other hand, the BilevelModel will require a solution \texttt{mode} which will select the technique used in the solution process. The final pieces of the basic interface are the \texttt{Lower} and \texttt{Upper} methods that direct JuMP macros to the proper bilevel optimization levels.

We will exemplify the basic interface by modeling the following simple bilevel optimization problem from \cite{dempe2002foundations}, Chapter 3.2, Page 25:
\[
\begin{aligned}
    &\min_{y \in {\mathbb{R}}} && 3x + y \\
    &\st && x \leq 5 \\
    &    && y \leq 8 \\
    &    && y \geq 0 \\
        &&& x(y) \in
     \begin{aligned}[t]
        &\argmin_{x \in {\mathbb{R}}} && -x\\
            &\st && x + y \leq 8\\
            &    && 4x + y \geq 8\\
            &    && 2x + y \leq 13\\
            &    && 2x - 7y \leq 0\\
     \end{aligned}\\
\end{aligned}
\]
The code to model, solve and query solutions is presented in Figure \ref{fig:jump_example}.
\begin{figure*}[!ht]
    \centering
\begin{lstlisting}[language = Julia]
using JuMP, BilevelJuMP, SCIP
model = BilevelModel(SCIP.Optimizer, mode = BilevelJuMP.SOS1Mode())
@variable(Upper(model), y)
@variable(Lower(model), x)
@objective(Upper(model), Min, 3x + y)
@constraints(Upper(model), begin
    x <= 5
    y <= 8
    y >= 0
end)
@objective(Lower(model), Min, -x)
@constraints(Lower(model), begin
     x +  y <= 8
    4x +  y >= 8
    2x +  y <= 13
    2x - 7y <= 0
end)
optimize!(model)
objective_value(model) # = 3 * (3.5 * 8/15) + 8/15
value(x) # = 3.5 * 8/15
value(y) # = 8/15
\end{lstlisting}
    \caption{Code to solve the example of \cite{dempe2002foundations}, Chapter 3.2, Page 25}
    \label{fig:jump_example}
\end{figure*}

We can follow the general workflow: include packages; initialize the model jointly with a solver, \texttt{SCIP} \citep{BestuzhevaEtal2021OO} in this case, and the solution mode \texttt{SOS1Mode} (modes will be discussed in the following sections); add variables to the proper levels so that they can be used by all constraints and objectives; add constraints and objectives to the proper levels (which can be done in any order); optimize the model; and query solutions.

\subsection{Solving Bilevel Optimization with Reformulations}

MathOptInterface defines a unique and well-posed interface that makes it possible to perform reformulations in problem instances to convert from one format into others. We start from an arbitrary user formulation in JuMP stored as an MOI model, then we can rewrite this model in an alternate form, which will lead to an MOI model, and consequently, we can pass the model to a solver wrapper that implements the MOI API. The solver optimizes the model and returns the solutions, which flows back to JuMP by applying the necessary mappings and transformations. The simplest and most used of these transformations are bridges \citep{legat2020mathoptinterface}, which are applied in individual variables, constraints, and objectives. The bridge system automatically converts a problem in the specific form expected by the solver.

However, some transformations require looking at the model as a whole and not only at its pieces (variables, objectives, and constraints). The first implementation of a whole model transformation was Dualization.jl \citep{guilherme_bodin_2021_4718987}. Dualization.jl's main function receives an MOI model and writes its dual in a second MOI model. Clearly, to perform such modification, the complete model must be known in advance. This feature is especially useful because some solvers only accept specific forms of mathematical programs. Hence, we can convert between primal and dual forms and solve the converted form without relying on the bridge system that might increase the problem size to reach the solver required form.

Dualization.jl also plays a key role in BilevelJuMP.jl's reformulations. Given a primal model like \eqref{ep:primal}, where $x$ are variables and $z$ are parameters, we can obtain the dual form following \cite{vandenberghe2010cvxopt}:
\begin{align}
& \max_{y_1, \ldots, y_m, w} & - \frac{1}{2} w^{\top} Q_1 w + \frac{1}{2} z^{\top} Q_3 z -\sum_{i=1}^m (b_i + D_i z)^{\top} y_i + d_0^{\top} z + b_0 \notag
\\
& \;\;\text{s.t.} & a_0 + Q_2 z - \sum_{i=1}^m A_i^{\top} y_i + Q_1 w & = 0
\label{eq:dual_feas}\\
& & y_i & \in \mathcal{C}_i^* & i = 1 \ldots m \label{eq:dual_cone}
\end{align}
We observe that the Dual feasibility constraints \eqref{eq:dual_feas}--\eqref{eq:dual_cone} are structurally very similar to the Dual Feasibility \eqref{eq:kkt_dual} and Stationarity \eqref{eq:kkt_stat} constraint sets from the KKT conditions. The only difference is that, in \eqref{eq:dual_feas}, $Q_1$ multiplies an additional variable $w$, and  in \eqref{eq:kkt_stat} $Q_1$ multiplies a parameter $z$.

BilevelJuMP.jl performs a more complicated model transformation. The two JuMP models that described each level of the bilevel program must be combined in a particular way to create the corresponding MPEC. Even though each variable belongs to one level, they are created in both but tagged with additional data to mark their level and their corresponding variable in the other model, constraints, and objectives. However, they only exist in the level they were created.

The first part of the transformation is to copy the upper level into a new model to append the other pieces of the MPEC. The second step is to add the KKT conditions of the second level. The Primal Feasibility constraints of the lower level are added as new constraints to the model (using the variable map between the two original JuMP models). Then the lower level model is dualized, considering upper-level variables as constants, and its constraints are passed to the new model to represent Stationarity and Dual Feasibility constraints. The additional variables, $w$, created in the dual problem are mapped into the upper variables $x$. At this point, we only need to add complementarity constraints.

\section{KKT Formulations}

In this section, we describe the reformulation of the conic MPEC to obtain mathematical programs that can be passed to existing solvers.

\subsection{Complementary slackness reformulations}

Starting from a convex problem, all the KKT conditions lead to convex constraints except the complementary slackness constraints. The main challenge in KKT reformulations is dealing with such non-linearity. Now we present some possible formulations which were already implemented and tested in BilevelJuMP.

We will assume that $y_i$ and $A_i x + b_i + D_i z$ are scalars, since almost all formulations rely on this assumption. The following formulations are restricted to: $\mathcal{C}_i \in \{ \mathbb{R}_+, \mathbb{R}_-, \{0\}  \}$ . Without loss of generality we will assume $\mathcal{C}_i = \mathbb{R}_+$.

\subsubsection{Special Ordered Sets of type 1}\label{sec-sos}
\

One SOS1 reformulation was presented in \cite{kleinert2020why} and in \cite{siddiqui2013sos1}. In BilevelJuMP.jl this formulation consists in replacing the complementarity constraints by the following:
\begin{align}
    & f_i = A_i x + b_i + D_i z \\
    & [y_i; f_i] \in SOS1
\end{align}

Considering this is the classic SOS1 set from \cite{beale1970special}, the SOS1 constraint implies that a solution is feasible only if at most one of the variables in the SOS1 set is different from zero.
It is equivalent to the original formula because one of the two scalars will have to be zero to have the product equal to zero.

Many solvers can handle this kind of constraint, e.g., Cbc, CPLEX, Gurobi, SCIP, Xpress, which makes this formulation particularly useful for practitioners.

\subsubsection{Indicator constraints}\label{sec-ind}
\

Indicator constraints \citep{belotti2016handling} and SOS1 are deeply related. A typical indicator constraint is defined by:
\begin{align}
    & x = 0 \implies Ay \leq b
\end{align}
This means that the constraint $Ay \leq b$ is only considered if $x = 0$, where $x$ is a binary variable. Analogously, another indicator constraint could depend on $x = 1$. Hence, one possible formulation for the complementarity slackness with indicator constraints is:
\begin{align}
    & f = 0 \implies A_i x + b_i + D_i z = 0 \\
    & f = 1 \implies y_i = 0 \\
    & f \in \{0,1\}
\end{align}

Many solvers are also capable of handling this kind of constraint which also makes this formulation very useful. As a final note on this formulation we note that a solver might not support Indicator Constraints for both $f = 0$ and $f = 1$, in this case we simply need one additional variable $g$ and the constraint: $f + g = 1$.

\subsubsection{Fortuny-Amat and McCarl}\label{sec-fa}
\

This formulation is commonly known by the name of the authors of \cite{fortuny1981representation} and is extremely used in practice. In very few words, it is an application of the big-M method:
\begin{align}
    & A_i x + b_i + D_i z \leq M_p f \\
    & y_i \leq M_d (1 - f) \\
    & f \in \{0,1\}
\end{align}
In such formulation, $M_p$ and $M_d$ are large numbers. We have assumed that both $A_i x + b_i + D_i z$ and $y_i$ are positive, thus, for each value of $f$ one of the elements in the complementarity pair is forced to zero.

The main drawback of this method is that the values of $M_p$ and $M_d$ must be large enough so that the optimal solution of the problem is not excluded. One can usually develop bounds on primal variables because the variable might be bounded due to the problem physics. However, finding reasonable bounds for dual variables might be much harder on specific applications. The work by \cite{pineda2019solving} shows that commonly used heuristics to select the big-Ms can fail. \cite{kleinert2020there} go further and demonstrate that verifying big-Ms is at least as hard as solving the Bilevel Problem itself.

When good bounds are available, the Fourtuny-Amat and McCarl formulation is very efficient in practice \citep{kleinert2020why}. Moreover, no extra constraints are required from solvers. Therefore, less complete MIP solvers like GLPK can be used to solve this kind of reformulation. On the other hand, the difficulty of computing bounds makes the SOS1 and Indicator formulations very useful for experimentation.

\subsubsection{Products}\label{sec-prod}
\

This is not a reformulation, because, in this case, the actual complementarity constraint in its product form is added to the optimization problem:
\begin{align}
& y_i^{\top} (A_i x + b_i + D_i z) = 0 \label{eq:prod}
\end{align}
NLP solvers frequently use this form to seek local optimal solutions. Although no guarantees of global optimality are provided when using this method, it is useful the get initial solutions to be used as bounds or even for cases where MIP solving is not practical. An additional weakness of this method is that \eqref{eq:prod} does not satisfy constraint qualification \citep{scholtes2001convergence, ralph2004some, fletcher2004solving} and is regularised as:
\begin{align}
& y_i^{\top} (A_i x + b_i + D_i z) \leq t
\end{align}
where $t$ is a small number.

In theory, one could reformulate all the products with binary expansion techniques such as the one in \cite{andrade2019enhancing} and use MIP solvers jointly with NLP solvers to reach solutions close to global optimal solutions. In practice, binary expansions also require bounds on the variables that are multiplied. This adds complications to the solution method because these cannot be added as regular constraints in the lower level; otherwise they would be dualized leading to more unbounded variables in both sides.

The binary expansion technique was implemented in QuadraticToBinary.jl  \citep{joaquim_dias_garcia_2021_4718981}, which can be used as an intermediary layer between BilevelJuMP.jl (or JuMP) and the selected solver—allowing any MIP solver with an MOI interface to solve approximations of quadratically constrained problems.

This formulation easily extends to vector sets. Hence, conic bilevel problems will require this formulation in BilevelJuMP.

\subsubsection{Complements}\label{sec-comp}
\

Some solvers are able to handle explicit complement constraints like Knitro \citep{nocedal2006knitro}, filterMPEC \citep{fletcher2021filterMPEC, fletcher2004solving}, NLPEC \citep{gams2021nlpec}. These solvers receive the constraints as special structures: pairs of variables or variable-expression pairs. Internally, the solver will employ their own reformulations.
\begin{align}
& y_i \perp A_i x + b_i + D_i z \label{eq:perp}
\end{align}

\subsubsection{Mixed mode}\label{sec-mixed}
\

Usually, practitioners select one single formulation and apply them to all complementarity constraints in the problem, but this is not a technical requirement. Consequently, one could combine formulations and select which formulation will be used for each pair. For instance, if one has good bounds for a specific pair, then just use Fortuny-Amat and McCarl for that constraint, while the other constraints would be reformulated with SOS1, for instance. Alternatively, even a conic bilevel with multiple linear constraints could be reformulated with SOS1 for all linear constraints and product mode (and binary expansions) for the conic constraints. We present an application of this method in Appendix \ref{ap-conic-mixed}.

\subsection{Primal Dual Equality reformulation}\label{sec-sd}

This formulation takes advantage of the fact that, under strong duality, the complementarity constraints are equivalent to enforcing that the primal and dual objective values are the same for a solution that is both primal and dual feasible. Therefore, the complementarity constraints are replaced by:
\begin{align}
     \frac{1}{2} x^{\top} P_1 x + x^{\top} P_2 z  + a_0^{\top} x
     =
     - \frac{1}{2} w^{\top} P_1 w  -\sum_{i=1}^m (b_i + D_i z)^{\top} y_i
\end{align}
Where the identical terms were already eliminated. This is also a non-convex quadratic constraint, even if the problem is linear due to $z$ and $y$ products.

One exciting feature of this formulation is that all the quadratic terms are concentrated in a single constraint, and the number of variable products might be much smaller than the number of complementarity constraints. Consequently, binary expansions were shown to be helpful to replace the quadratic terms and achieve approximate global optimal solutions in \cite{pereira2005strategic, zare2019note}.

\subsection{Comparison of methods}

We present a brief comparison between the solution methods. Table \ref{tab-reform} presents the method name in BilevelJuMP, section in which it is described, solver requirements and additional comments.

\begin{table}[h]
\resizebox{\textwidth}{!}{
\begin{tabular}{l|l|l|l}
Method Name       & Sec. & Solver requirement                                         & Comments \\
\hline
\texttt{SOS1Mode}          & \ref{sec-sos}       & MIP solver with SOS of type 1                              & No additional information. Only linear constraints. \\
\texttt{IndicatorMode}     & \ref{sec-ind}       & MIP solver with Indicator Constraints                      & No additional information. Only linear constraints. \\
\texttt{FortunyAmatMcCarlMode}& \ref{sec-fa}      & MIP solver                                                 & Require non-trivial big-M. Only linear constraints. \\
\texttt{ProductMode}       & \ref{sec-prod}      & Non-convex quadratic constraints                           & Works with conic constraints. Require regularization.\\
\texttt{ComplementMode}    & \ref{sec-comp}   & Complementarity constraint                             & Few solvers supporting such constraints.\\
\texttt{MixedMode}         & \ref{sec-mixed}     & Requirements of selected methods & Pros and cons from selected methods.\\
\texttt{StrongDualityMode} & \ref{sec-sd}      & Non-convex quadratic constraints                           & Works with conic constraints.
\end{tabular}
}
\caption{Reformulation methods}
\label{tab-reform}
\end{table}

\section{Example}

In this section, we describe a slightly more interesting example of a bilevel program. The main goal is to start from a non-trivial problem, model it in BilevelJuMP, and solve it with multiple methods to have a glimpse of the multitude of applications of the package.

As previously mentioned, hyperparameter tuning with bilevel optimization is a recent trend in the intersection of the Machine Learning and Optimization communities \citep{franceschi2018bilevel, kunisch2013bilevel, mackay2019self}. Although most hyperparameter tuning methods based on bilevel optimization are usually heuristic with special considerations to the problem in question, this is a good case to describe the functionality of the package due to the simplicity of the model and because small enough instances can be solved by standard methods implemented in this package.

We have selected hyperparameter tuning in support vector regressions (SVR). The example will follow the one from \cite{bennett2006model}, though with some simplifications. Given two data sets $O$ and $I$ with out-of-sample and in-sample data represented by the points labelled by $i$: $(y_i, \{x_{ij}\}_{j \in J})$, $J$ is the feature set.
\[
\begin{aligned}
    &\min_{C \geq 0, \varepsilon \geq 0, \xi^U \geq 0} && \sum_{i \in O} \xi_i^U \\
    &\st && \xi_i^U \geq  + y_i - \sum_j w_j x_{ij}, \quad i \in O \\
    &    && \xi_i^U \geq - y_i + \sum_j w_j x_{ij} , \quad i \in O \\
        &&& w(C, \varepsilon) \in
     \begin{aligned}[t]
        &\argmin_{\xi^L \geq 0, w } && ||w||^2_2 + C \sum_{i \in I} \xi_i^L\\
            &\st && \xi_i^L + \varepsilon \geq  + y_i - \sum_{j \in J} w_j x_{ij} , \quad i \in I\\
            &    && \xi_i^L + \varepsilon \geq - y_i + \sum_{j \in J} w_j x_{ij}, \quad i \in I\\
     \end{aligned}\\
\end{aligned}
\]
The lower model is responsible for obtaining the best possible support vectors $w$ given the problem data and the hyperparameters $C$ and $\varepsilon$, which are variables selected by the upper level so that the $w$ optimized by the lower level has minimal out-of-sample error.
The variables $\xi^U$ and $\xi^L$ denote the absolute value loss in the upper and lower models, respectively. The upper level is a linear program, while the lower level is quadratic. In Figure \ref{fig:jump_svr} we present BilevelJuMP.jl code to model the hyperparameter tuning of SVR described above. Thanks to the JuMP syntax, the code greatly resembles the abstract model, simplifying the writing and documenting of the code. 
\begin{figure*}[!ht]
    \centering
\begin{lstlisting}[language = Julia]
using JuMP, BilevelJuMP
# sample data
Features = 2
Samples = 10
J = 1:Features
I = 1:div(Samples, 2)
O = (div(Samples, 2)+1):Samples
x = 2 * (rand(Samples, Features) .- 0.5)
w_real = ones(Features)
y = x * w_real .+ 0.1 * 2 * (rand(Samples) .- 0.5)
# model building
model = BilevelModel()
@variable(Upper(model), C >= 0)
@variable(Upper(model), eps >= 0)
@variable(Upper(model), xi_U[i=O] >= 0)
@variable(Lower(model), w[j=J])
@variable(Lower(model), xi_L[i=I] >= 0)
@objective(Upper(model),
    Min, sum(xi_U[i] for i in O))
@constraints(Upper(model), begin
    [i in O], xi_U[i] >= + y[i] - sum(w[j]*x[i,j] for j in J)
    [i in O], xi_U[i] >= - y[i] + sum(w[j]*x[i,j] for j in J)
end)
@objective(Lower(model),
    Min, sum(w[j]^2 for j in J) + C * sum(xi_L[i] for i in I))
@constraints(Lower(model), begin
    [i in I], xi_L[i] + eps >= + y[i] - sum(w[j]*x[i,j] for j in J)
    [i in I], xi_L[i] + eps >= - y[i] + sum(w[j]*x[i,j] for j in J)
end)
\end{lstlisting}
    \caption{Code to model SVR hyper-parameter tuning}
    \label{fig:jump_svr}
\end{figure*}

That same code was used to perform a series of comparisons between solvers. We started by creating instances with a different number of features and observations (dataset size). We randomly created the matrix $x$ with a uniform distribution in $[-1, +1]$, then we created the \textit{real} $w$ as a vector of ones with appropriate dimension. Next, we defined $y = x w + \epsilon$, where $\epsilon$ follows a uniform distribution in $[-0.1, +0.1]$. Half of the dataset was considered in-sample data, while the other half was considered out-of-sample data. It is not our intention to be fully realistic here, our goal is to provide a didatic example.

We created instances with $10$, $100$ and $1000$ samples. For all these sample sizes we created samples with $1$, $2$ and $5$ features. For the datasets with $100$ samples we also created datasets with $10$, $20$ and $50$ features.

Finally, we optimized the bilevel problem for each data set with multiple reformulations and with multiple solvers. The only solver attribute we set was a time limit of $600$ seconds ($10$ minutes) and left all other attributes as default, which might differ considerably from one solver to the other. Again, our primary goal is not a detailed and rigorous comparison of solvers but to show the software's functionality in a usage example as practitioners and researchers might want to solve the same problem with multiple methods and select the one that best fits their needs.

We present results in the following tables. We used
Julia 1.6.2,
CPLEX 22.1 \citep{cplex},
Gurobi 9.5 \citep{Gurobi},
HiGHS 1.2 \citep{huangfu_2018},
Ipopt 3.14 \citep{wachter2006implementation},
Knitro 13.0 \citep{nocedal2006knitro},
SCIP 8.0 \citep{BestuzhevaEtal2021OO},
Xpress 8.13 \citep{xpress}. All the required code is in the benchmarks folder of the git repository, including exact package versions (see the \textit{manifest.toml} file).

All tables have a similar format. The first column describes the instance, the first number being the sample size and the second the number of features. Then we have three columns for each solver, the first, \textit{Obj}, presents the upper-level objective value returned by the solver (typically the best incumbent solution), the second contains the \textit{Gap} in percent (\%), Ipopt and KNITRO will not have gaps as they are NLP solvers, if no gap was reported the entry will be blank (with a ``$-$''), the third is \textit{Time} in seconds, if the time reaches $600$ the entry will be blank (with a ``$-$'').

Table \ref{table_sos1} presents results for  \texttt{SOS1Mode} and \texttt{IndicatorMode}. Table \ref{table_fa100} presents results for \texttt{FortunyAmatMcCarlMode}, with big-Ms set to $100$, \texttt{StrongDualityMode} and \texttt{ProductMode}, the latter two with binary expansions, so the resulting problem is a MIP, where the variable bounds were set to  $+/-100$. Finally, Table \ref{table_prod} presents the solutions of both \texttt{ProductMode} and \texttt{StrongDualityMode} for Non-Linear Programming solvers and Gurobi with its \textit{NonConvex} mode activated.

We can draw some conclusions from the tables. We note that \texttt{SOS1Mode} and \texttt{IndicatorMode} perform well in smaller instances, with a slight advantage for \texttt{SOS1Mode}. Interestingly, CPLEX's solution for $1000/01$ with \texttt{IndicatorMode} slightly disagrees with all solvers solution with the \texttt{SOS1Mode}. \texttt{FortunyAmatMcCarlMode} and \texttt{StrongDualityMode} seem very amenable to MIP solvers with Gurobi closing the gap within the given 10 minutes for all but one instance in the latter mode. However, we must be careful since we selected arbitrary bounds for those methods and \texttt{StrongDualityMode} also relies on binary expansion approximations, which led solvers to a solution that disagrees with the other methods on the $10/05$ instance. On the other hand, \texttt{ProductMode} is the worst strategy for MIP solvers in these instances. For NLP solvers both \texttt{ProductMode} and \texttt{StrongDualityMode} return objective values that are close to the ones found by MIP solvers, but in this case there is a slight advantage for \texttt{ProductMode}. Finally, Gurobi NonConvex seems to work much better with \texttt{StrongDualityMode} claiming very good results in the instances with $1000$ samples that agree with some of the other presented objective values.

The results are particular to a toy problem. However, the tables demonstrate that the software can interface with multiple solvers and consider multiple methods. Moreover, there is value in experimenting with multiple solvers and methods implemented in BilevelJuMP.

\input{table2_sos1}

\input{table2_fa100}

\input{table2_prod}

\clearpage

\newpage

\section{Package Comparison}

A comparison between BilevelJuMP.jl and four other bilevel optimizations modeling interfaces that include solution methods is presented in Table \ref{table_software}. We include BilevelOptimization.jl as it was the key motivation for BilevelJuMP.jl; PAO, as the new bilevel interface of pyomo; GAMS that relies on EMP; and YALMIP that motivated the development of Dualization.jl.

In each table line, we briefly depict the answer to each of the following questions:

\begin{enumerate}
  \item Which programming language does a user have to code the models?
  \item What is the licensing scheme? (MIT is the most permissive among the ones shown).
  \item Does the modeling interface support MIP solver-based methods, like SOS1 and big-M?
  \item Does the modeling interface support NLP solver-based methods, like products?
  \item Does the modeling interface support MPEC solvers that accept explicit complementarity constraints?
  \item Can the user access dual variables of the lower-level problem and use them explicitly while modeling the upper-level problem?
  \item Which problem classes are accepted in the lower-level problem?
  \item Which problem classes are accepted in the upper-level problem?
\end{enumerate}

The two latter questions used the following code: CP is Conic Programming, QP is linear programming with optional quadratic objective, NLP stands for Non-Linear Programming, VI represents Variational Inequalities, and Int is Integer Programming. Although one can model problems of given classes, specific solvers will be required to handle the resulting formulations. Finally, we note that all classes might not be supported simultaneously by all the interfaces, BilevelJuMP.jl supports all the described classes in the same model. Finally, it is worth mentioning that the possibility to consider bilevel models in which lower-level primal and dual variables are present in the first-level problem significantly enlarges the spectrum of practical applications that can be covered with the package. We marked YALMIP as ready to handle lower level duals because this can be achieved by explicitly calling their ``kkt'' function on the lower level data and appending to the primal problem. For instance, strategic bidding as well as market-power assessments in electricity markets highly depend on bilevel models with such dependencies \cite{fanzeresahmed2019}.

\begin{table}[]
\centering
\resizebox{16cm}{!}{
\begin{tabular}{l|lllll}
Name         & BilevelJuMP.jl   & BilevelOptimization.jl & PAO/Pyomo & GAMS       & YALMIP        \\ \hline
Language     & Julia         & Julia                  & Python    & GAMS       & MATLAB        \\
License      & MIT           & MIT                    & BSD       & Commercial & YALMIP        \\
MIP solvers  & Yes           & Yes                    & Yes       & No         & Yes           \\
NLP solvers  & Yes           & No                     & Yes       & No         & Yes           \\
MPEC solvers & Yes           & No                     & No        & Yes        & Yes           \\
DualVar      & Yes           & No                     & No        & Yes        & Yes           \\
Lower Level  & CP/QP         & QP (matrix form)       & QP/Int    & QP/NLP/VI  & QP            \\
Upper Level  & CP/QP/NLP/Int & CP/QP/NLP/Int          & QP/Int    & QP/NLP/Int & CP/QP/NLP/Int
\end{tabular}
}
\caption{Modeling interfaces for bilevel optimization} \label{table_software} 
\end{table}

\section{Conclusion}\label{conclusion}

We presented BilevelJuMP, an open-source Julia library for Bilevel Optimization that allows the user to model a wide range of bilevel optimization problems very easily. Moreover, the user has access to multiple reformulation techniques that can be considered to handle different problems better. More specifically, BilevelJuMP.jl allows modeling very general problems in the upper level (all JuMP supported formulations, such as non-linear, conic, mixed integer constraints) and conic problems in the lower level. Additionally, it implements multiple MPEC based reformulation techniques and MIP or NLP key as solution algorithms. This broad and flexible infrastructure of models and methods all built in a single open- source package for JuMP allows practitioners to use BilevelJuMP.jl for quick deploy and run experiments using state-of-the-art solvers and methods. It can be used by students introduced to the realm of bilevel optimization due to its easy-to-use and flexible structures, researchers and developers that can quickly test (or benchmark) new methods and models, and also develop new applications and packages, as well as by industry practitioners, who may not be familiar with bilevel solution strategies, but can rely on the package to address specific bilevel problems composing parts of real-world applications.

BilevelJuMP.jl joins a group of JuMP and MOI extensions that were made possible thanks to the good design of the latter two.

Just like JuMP and MOI, BilevelJuMP.jl is under active development, and more features are planned to be included. The library has gotten great interest from other contributors that are currently working towards new features, including support for integer variables in the lower-level with the solver MibS \citep{tahernejad2020branch}. Further developments include: implementing other techniques such as valid inequalities \citep{kleinert2021closing}, column-constraint generation based techniques \citep{yue2019projection}; developing a file format for bilevel optimization based on MathOptFormat; integrating other solvers such as the one from \cite{fischetti2020solver}.

\ 

\ACKNOWLEDGMENT{%
The authors thank the JuMP and MOI contributors for making such a great tool available. The authors also thank all the contributors and early users of BilevelJuMP, especially Mathieu Besaçon, Hesam Shaelaie, and Oscar Dowson. Authors were partially supported by the Coordenação de Aperfeiçoamento de
Pessoal de Nível Superior - Brasil (CAPES) - Finance Code 001. The work of Alexandre Street was also partially supported by FAPERJ and CNPq.
}

%
%
%

\

\newpage

\

\begin{APPENDICES}
\section{Lower-level duals}\label{ap-dual}

This modeling feature enables the implementation of workflows where one (or more) of the upper-level variables is the dual of a lower-level constraint. In particular, in the energy sector, it is common to model the energy prices as the dual variable associated with the energy demand equilibrium constraint. One example of an application that uses this feature is \cite{fanzeresahmed2019}, which focuses on strategic bidding in auction-based energy markets.
A small and simplified example of the modeled problem would be the model:
\begin{align}
    &\max_{\lambda, q_S} \quad \lambda \cdot g_S \\
    &\st \quad 0 \leq q_S \leq 100\\
    &\hspace{28pt} (g_S, \lambda) \in \argmin_{g_S, g_{1}, g_{2}, g_D} 50 g_{R1} + 100  g_{R2} + 1000 g_{D}\\
            & \hspace{70pt} \st \quad g_S \leq q_S \\
            & \hspace{88pt} \quad  0 \leq g_S \leq 100 \\
            & \hspace{88pt}\quad  0 \leq g_{1} \leq 40 \\
            & \hspace{88pt}\quad  0 \leq g_{2} \leq 40 \\
            & \hspace{88pt}\quad  0 \leq g_{D} \leq 100 \\
    & \hspace{88pt}\quad  g_S + g_{1} + g_{2} + g_D = 100 \quad  : \quad \lambda \label{eq-dual-lambda}
\end{align}
Where $\lambda$ is the dual of the load balance constraint \eqref{eq-dual-lambda}, $g_S$, $g_{1}$, $g_2$ represent the generation of the strategic bidder and from two other (non-strategic) plants. $g_D$ represents the deficit in generation. Finally, $q_S$ is the quantity bid optimized by the strategic generator.

BilevelJuMP.jl allows users to implement similar models using the function \texttt{DualOf} that binds a new variable in the upper level to an existing constraint in the lower level.
The model can be written as:

\begin{lstlisting}[language = Julia]
model = BilevelModel()
@variable(Upper(model), 0 <= qS <= 100)
@variable(Lower(model), 0 <= gS <= 100)
@variable(Lower(model), 0 <= gR1 <= 40)
@variable(Lower(model), 0 <= gR2 <= 40)
@variable(Lower(model), 0 <= gD <= 100)
@objective(Lower(model), Min, 50gR1 + 100gR2 + 1000gD)
@constraint(Lower(model), gS <= qS)
@constraint(Lower(model), demand_equilibrium, gS + gR1 + gR2 + gD == 100)
@variable(Upper(model), lambda, DualOf(demand_equilibrium))
@objective(Upper(model), Max, lambda*gS)
\end{lstlisting}

\subsection{NLP solution}

This model, can the be solved by selecting a reformulation and a solver. Here we select Strong-Duality reformulation, the Ipopt solver and call optimizes to perform the reformulation and solve it.

\begin{lstlisting}[language = Julia]
BilevelJuMP.set_mode(model, BilevelJuMP.StrongDualityMode())
set_optimizer(model, Ipopt.Optimizer)
optimize!(model)
\end{lstlisting}

\subsection{MIP solution}

It is also possible to solve such problem by using a MIP formulation. The main issue is the product of variable in the upper level objective. However, this can be easily handled by using the aforementioned \texttt{QuadraticToBinary} package for automatic binary expansions. Because binary expansions require bounds on variables, we change the following line:
\begin{lstlisting}[language = Julia]
@variable(Upper(model), 0 <= lambda <= 1000, DualOf(demand_equilibrium))
\end{lstlisting}
Then, as before, we set a solver (now SCIP with the \texttt{QuadraticToBinary} wrapper) and a solution method (now Fortuny-Amat and McCarl):
\begin{lstlisting}[language = Julia]
set_optimizer(model,
    ()->QuadraticToBinary.Optimizer{Float64}(SCIP.Optimizer()))
BilevelJuMP.set_mode(model,
    BilevelJuMP.FortunyAmatMcCarlMode(dual_big_M = 100))
optimize!(model)
\end{lstlisting}

\ 
\section{Conic Bilevel and Mixed Mode}\label{ap-conic-mixed}

Here we present a simple bilevel program with a conic lower level model described in example 3.3 from \cite{chi2014models}. 

\begin{align}
    &\max_{x \in \mathbb{R}} \quad x + 3y_1 \\
    &\st \quad 2 \leq x \leq 6\\
    & \hspace{28pt} y(x) \in \argmin_{y\in {\mathbb{R}^3}} -y_1\\
            & \hspace{58pt} \st \quad x + y_1 \leq 8 \\
            & \hspace{76pt} \quad x + 4y_1 \geq 8 \\
            & \hspace{76pt}  \quad x + 2y_1 \leq 12 \\
            & \hspace{76pt}  \quad y \in {SOC}_3 \label{eq-soc}
\end{align}

In this problem most of the constraints are regular linear constraints while the last one, \eqref{eq-soc}, is a second order cone constraint. Such constraint ensures that the vector $y$ belongs to a second order cone of dimension $3$, that is: $y_1 \geq \sqrt{y_2^2 + y_3^2}$. This problem can be encoded using regular JuMP syntax for conic programs:

\begin{lstlisting}[language = Julia]
model = BilevelModel()
@variable(Upper(model), x)
@variable(Lower(model), y[i=1:3])
@objective(Upper(model), Min, x + 3y[1])
@constraint(Upper(model), x >= 2)
@constraint(Upper(model), x <= 6)
@objective(Lower(model), Min, - y[1])
@constraint(Lower(model), con1, x +  y[1] <=  8)
@constraint(Lower(model), con2, x + 4y[1] >=  8)
@constraint(Lower(model), con3, x + 2y[1] <= 12)
@constraint(Lower(model), con4, y in SecondOrderCone())
\end{lstlisting}

\subsection{NLP solution and start values}

We can set, for instance, the product reformulation and selected Ipopt as a solver. As Ipopt does not have native support for second order cones, we use the non-default MOI bridge \texttt{SOCtoNonConvexQuad} to convert second order cones into quadratic constraints.

\begin{lstlisting}[language = Julia]
BilevelJuMP.set_mode(model,BilevelJuMP.ProductMode(1e-5))
set_optimizer(model,
    ()->MOI.Bridges.Constraint.SOCtoNonConvexQuad{Float64}(Ipopt.Optimizer()))
optimize!(model)
\end{lstlisting}

This problem is very simple, but more complex models might require more information such as starting points that can be passed on variable creation with standard JuMP syntax, for instance:
\begin{lstlisting}[language = Julia]
@variable(Upper(model), x, start = 6)
\end{lstlisting}
The user could also use the alternative JuMP syntax:
\begin{lstlisting}[language = Julia]
set_start_value(x, 6)
set_dual_start_value(con2, 0)
\end{lstlisting}

\subsection{MIP solution and mixed mode}

Alternatively, we could have used a Mixed Integer Second Order Cone Program (MISOCP) solver together with binary expansions. Complementarity of conic constraints is more difficult to handle because they require a sum of products that cannot be reformulated with other methods. Therefore, we rely on product reformulation for conic constraints. However, we can use other reformulations like indicator constraints for the non-conic constraints. Mixing the two of them can be done with Mixed Mode from Section \ref{sec-mixed}.

The following code describes how to solve the problem with a MISOCP based solver.

\begin{lstlisting}[language = Julia]
set_optimizer(model,
    ()->QuadraticToBinary.Optimizer{Float64}(Xpress.Optimizer(),lb=-10,ub=10))
BilevelJuMP.set_mode(model, 
    BilevelJuMP.MixedMode(default = BilevelJuMP.IndicatorMode()))
BilevelJuMP.set_mode(con4, BilevelJuMP.ProductMode(1e-5))
optimize!(model)
\end{lstlisting}

We set the reformulation method as Mixed Mode and selected Indicator constraints to be the default for the case in which we do not explicitly specify the reformulation. Then we set product mode for the second order cone reformulation.

As described in Appendix \ref{ap-dual}, binary expansions require bounded variables, hence the \texttt{QuadraticToBinary} meta-solver accepts fallback to upper and lower bounds (\texttt{ub} and \texttt{lb}), used for variables with no explicit bounds.

\end{APPENDICES}


\bibliographystyle{informs2014} 
\bibliography{ref.bib} 



\end{document}

%% file: table2_sos1.tex
\begin{table}[!ht]
\centering
\resizebox{13cm}{!}{
\begin{tabular}{rr|rrr|rrr|rrr|rrr}
\toprule
 &  & \multicolumn{3}{c}{CPLEX} \vline & \multicolumn{3}{c}{Gurobi} \vline & \multicolumn{3}{c}{SCIP} \vline & \multicolumn{3}{c}{Xpress}  \\
 & Inst & Obj & Gap  & Time  & Obj & Gap  & Time  & Obj & Gap  & Time  & Obj & Gap  & Time  \\
\midrule
 & 10/01  & 0.30 & 0 & 0 & 0.30 & 0 & 0 & 0.30 & 0 & 0 & 0.30 & 0 & 0 \\
 & 10/02  & 0.22 & 0 & 0 & 0.22 & 0 & 0 & 0.22 & 0 & 0 & 0.22 & 0 & 0 \\
 \parbox[t]{2mm}{\multirow{3}{*}{\rotatebox[origin=c]{90}{\texttt{SOS1Mode}}}}
 & 10/05  & 0.09 & 0 & 0 & 0.09 & 0 & 0 & 0.09 & 0 & 0 & 0.09 & 0 & 0 \\
 & 100/01  & 2.42 & 0 & 0 & 2.42 & 0 & 0 & 2.42 & 0 & 0 & 2.42 & 0 & 0 \\
 & 100/02  & 2.40 & 4 &  -  & 2.40 & 4 &  -  & 2.40 & 4 &  -  & 2.40 & 4 &  -  \\
 & 100/05  & 2.30 & 6 &  -  & 2.30 & 6 &  -  & 2.31 & 6 &  -  & 54.87 &  -  &  -  \\
 & 100/10  & 8.54 & 392 &  -  & 79.59 &  -  &  -  & 79.59 &  -  &  -  & 79.59 &  -  &  -  \\
 & 100/20  & 102.79 &  -  &  -  & 8.21 & 457 &  -  & 102.79 &  -  &  -  & 96.89 &  -  &  -  \\
 & 100/50  & 23.35 & 307 &  -  & 23.35 & 299 &  -  & 23.35 &  -  &  -  & 23.35 & 350 &  -  \\
 & 1000/01  & 25.02 & 0 &  -  & 28.63 & 14 &  -  & 25.02 & 0 &  -  & 25.02 & 0 &  -  \\
 & 1000/02  & 323.30 &  -  &  -  & 323.30 &  -  &  -  & 323.30 &  -  &  -  & 323.30 &  -  &  -  \\
 & 1000/05  & 533.37 &  -  &  -  & 533.37 &  -  &  -  & 533.37 &  -  &  -  & 533.37 &  -  &  -  \\
\midrule
 & 10/01  & 0.30 & 0 & 0 & 0.30 & 0 & 0 & 0.30 & 0 & 0 & 0.30 & 0 & 0 \\
 & 10/02  & 0.22 & 0 & 0 & 0.22 & 0 & 0 & 0.22 & 0 & 0 & 0.22 & 0 & 0 \\
 \parbox[t]{2mm}{\multirow{3}{*}{\rotatebox[origin=c]{90}{\texttt{IndicatorMode}}}}
 & 10/05  & 0.09 & 0 & 0 & 0.09 & 0 & 0 & 0.09 & 0 & 0 & 0.09 & 0 & 0 \\
 & 100/01  & 2.42 & 0 & 0 & 2.42 & 0 & 0 & 2.42 & 0 & 2 & 2.42 & 0 & 2 \\
 & 100/02  & 2.40 & 4 &  -  & 2.40 & 4 &  -  & 9.08 & 294 &  -  & 2.42 & 5 &  -  \\
 & 100/05  & 2.30 & 6 &  -  & 2.30 & 6 &  -  & 39.08 &  -  &  -  & 2.31 & 6 &  -  \\
& 100/10  & 79.59 &  -  &  -  & 79.59 &  -  &  -  & 79.59 &  -  &  -  & 79.59 &  -  &  -  \\
 & 100/20  & 102.79 &  -  &  -  &  -  &  -  &  -  & 102.79 &  -  &  -  & 102.79 &  -  &  -  \\
 & 100/50  & 23.35 &  -  &  -  & 23.35 &  -  &  -  & 23.35 &  -  &  -  & 23.35 & 576 &  -  \\
 & 1000/01  & 25.06 & 0 &  -  &  -  &  -  &  -  & 77.45 & 209 &  -  & 195.11 & 680 &  -  \\
 & 1000/02  & 323.30 &  -  &  -  &  -  &  -  &  -  & 323.30 &  -  &  -  &  -  &  -  &  -  \\
 & 1000/05  & 533.37 &  -  &  -  & 533.37 &  -  &  -  & 533.37 &  -  &  -  &  -  &  -  &  -  \\
\bottomrule
\end{tabular}
}
\caption{MIP solvers with \texttt{SOS1Mode} and \texttt{IndicatorMode}, Time in seconds (s), Gap in percent (\%).}
\label{table_sos1}
\end{table}

%% file: table2_fa100.tex
\begin{table}[!ht]
\resizebox{\textwidth}{!}{
\begin{tabular}{rr|rrr|rrr|rrr|rrr|rrr}
\toprule
 &  & \multicolumn{3}{c}{CPLEX} \vline & \multicolumn{3}{c}{Gurobi} \vline & \multicolumn{3}{c}{HiGHS} \vline & \multicolumn{3}{c}{SCIP} \vline & \multicolumn{3}{c}{Xpress}  \\
 & Inst & Obj & Gap  & Time  & Obj & Gap  & Time  & Obj & Gap  & Time  & Obj & Gap  & Time  & Obj & Gap  & Time  \\
\midrule
 & 10/01  & 0.30 & 0 & 0 & 0.30 & 0 & 0 & 0.30 & 0 & 0 & 0.30 & 0 & 0 & 0.30 & 0 & 0 \\
 & 10/02  & 0.22 & 0 & 0 & 0.22 & 0 & 0 & 0.22 & 0 & 0 & 0.22 & 0 & 0 & 0.22 & 0 & 0 \\
 \parbox[t]{2mm}{\multirow{3}{*}{\rotatebox[origin=c]{90}{\texttt{FortunyAmatMcCarlMode}}}}
 & 10/05  & 0.09 & 0 & 0 & 0.09 & 0 & 0 & 0.09 & 0 & 0 & 0.09 & 0 & 0 & 0.09 & 0 & 0 \\
 & 100/01  & 2.42 & 0 & 0 & 2.42 & 0 & 0 & 2.42 & 0 & 0 & 2.42 & 0 & 1 & 2.42 & 0 & 0 \\
 & 100/02  & 2.40 & 4 &  -  & 2.40 & 4 &  -  & 2.40 & 4 &  -  & 2.43 & 5 &  -  & 2.43 & 5 &  -  \\
 & 100/05  & 2.30 & 6 &  -  & 2.29 & 5 &  -  & 54.87 &  -  &  -  & 39.08 &  -  &  -  & 2.31 & 6 &  -  \\
 & 100/10  & 79.59 &  -  &  -  & 2.33 & 34 &  -  & 79.59 &  -  &  -  & 79.59 &  -  &  -  & 79.59 &  -  &  -  \\
 & 100/20  & 22.11 &  -  &  -  & 22.17 &  -  &  -  & 102.79 &  -  &  -  & 102.79 &  -  &  -  &  -  &  -  &  -  \\
 & 100/50  & 23.35 & 355 &  -  & 23.35 & 330 &  -  & 23.35 &  -  &  -  & 23.35 & 952 &  -  & 23.35 & 736 &  -  \\
 & 1000/01  & 25.02 & 0 &  -  & 25.02 & 0 &  -  & 25.02 & 0 &  -  & 70.13 & 180 &  -  & 25.02 & 0 &  -  \\
 & 1000/02  & 24.46 & 3 &  -  & 23.74 & 0 & 12 & 323.30 &  -  &  -  & 323.30 &  -  &  -  & 23.75 & 0 &  -  \\
 & 1000/05  & 533.37 &  -  &  -  & 533.37 &  -  &  -  & 533.37 &  -  &  -  & 533.37 &  -  &  -  & 533.37 &  -  &  -  \\
  \midrule
 & 10/01  & 0.30 & 0 & 344 & 0.30 & 0 & 17 & 0.30 & 0 & 269 & 0.30 & 0 & 83 & 0.30 & 0 & 447 \\
 & 10/02  & 1.82 & 739 &  -  & 0.22 & 0 &  -  & 3.12 &  -  &  -  & 0.33 & 0 & 495 & 0.30 & 38 &  -  \\
 \parbox[t]{2mm}{\multirow{3}{*}{\rotatebox[origin=c]{90}{\texttt{ProductMode}}}}
 & 10/05  & 0.53 &  -  &  -  & 7.22 &  -  &  -  & 8.72 &  -  &  -  & 0.67 &  -  &  -  & 0.60 &  -  &  -  \\
 & 100/01  & 18.08 & 647 &  -  & 14.37 & 494 &  -  & 23.28 & 0 & 52 & 21.87 & 803 &  -  & 20.74 & 0 & 55 \\
 & 100/02  & 27.06 &  -  &  -  &  -  &  -  &  -  &  -  &  -  &  -  &  -  &  -  &  -  &  -  &  -  &  -  \\
 & 100/05  &  -  &  -  &  -  & 48.56 &  -  &  -  &  -  &  -  &  -  &  -  &  -  &  -  &  -  &  -  &  -  \\
 & 100/10  &  -  &  -  &  -  & 75.06 &  -  &  -  &  -  &  -  &  -  &  -  &  -  &  -  & 78.62 &  -  &  -  \\
 & 100/20  &  -  &  -  &  -  & 99.82 &  -  &  -  &  -  &  -  &  -  &  -  &  -  &  -  & 101.51 &  -  &  -  \\
 & 100/50  & 247784.27 &  -  &  -  & 183.62 &  -  &  -  &  -  &  -  &  -  &  -  &  -  &  -  &  -  &  -  &  -  \\
 & 1000/01  & 45.15 & 80 &  -  &  -  &  -  &  -  &  -  &  -  &  -  & 58.77 & 135 &  -  & 147.68 & 490 &  -  \\
 & 1000/02  &  -  &  -  &  -  &  -  &  -  &  -  &  -  &  -  &  -  &  -  &  -  &  -  & 165.04 & 595 &  -  \\
 & 1000/05  &  -  &  -  &  -  &  -  &  -  &  -  &  -  &  -  &  -  &  -  &  -  &  -  &  -  &  -  &  -  \\
 \midrule
 & 10/01  & 0.30 & 0 & 512 & 0.30 & 0 &  -  & 0.30 & 0 &  -  & 0.30 & 0 &  -  & 0.30 & 0 &  -  \\
 & 10/02  & 0.22 & 0 & 0 & 0.22 & 0 & 0 & 0.22 & 0 & 7 & 0.22 & 0 & 36 & 0.22 & 0 & 52 \\
 \parbox[t]{2mm}{\multirow{3}{*}{\rotatebox[origin=c]{90}{\texttt{StrongDualityMode}}}}
 & 10/05  & 0.00 & 0 & 2 & 0.00 & 0 & 3 & 0.00 & 0 & 165 & 0.00 & 0 & 37 & 0.00 &  -  & 4 \\
 & 100/01  & 2.42 & 0 & 2 & 2.42 & 0 & 1 & 2.42 & 0 & 2 &  -  &  -  &  -  & 2.42 & 0 & 0 \\
 & 100/02  & 2.30 & 0 & 7 & 2.30 & 0 & 10 & 2.30 & 0 & 63 & 2.36 & 2 &  -  &  -  &  -  &  -  \\
 & 100/05  & 2.16 & 0 & 139 & 2.16 & 0 & 66 & 41.77 &  -  &  -  &  -  &  -  &  -  & 2.18 & 0 &  -  \\
 & 100/10  & 1.73 & 0 & 269 & 1.73 & 0 & 34 & 75.84 &  -  &  -  & 78.61 &  -  &  -  &  -  &  -  &  -  \\
 & 100/20  & 90.35 &  -  &  -  & 1.45 & 0 & 269 &  -  &  -  &  -  &  -  &  -  &  -  & 1.57 & 8 &  -  \\
 & 100/50  & 176.21 &  -  &  -  & 160.70 &  -  &  -  &  -  &  -  &  -  &  -  &  -  &  -  &  -  &  -  &  -  \\
 & 1000/01  & 25.01 & 0 & 153 & 25.01 & 0 & 22 & 25.01 & 0 & 32 & 25.01 & 0 & 20 & 25.01 & 0 & 10 \\
 & 1000/02  & 23.74 & 0 & 39 & 23.74 & 0 & 23 & 23.75 & 0 &  -  &  -  &  -  &  -  & 23.74 & 0 & 20 \\
 & 1000/05  & 529.71 &  -  &  -  & 24.43 & 0 & 543 & 408.46 &  -  &  -  &  -  &  -  &  -  & 24.43 & 0 & 482 \\
\bottomrule
\end{tabular}
}
\caption{MIP solvers with \texttt{FortunyAmatMcCarlMode}, \texttt{ProductMode} and \texttt{StrongDualityMode}, Time in seconds (s), Gap in percent (\%).}
\label{table_fa100}
\end{table}

%% file: table2_prod.tex
\begin{table}[!ht]
\centering
\resizebox{11cm}{!}{
\begin{tabular}{rr|rrr|rrr|rrr}
\toprule
 &  & \multicolumn{3}{c}{Gurobi NonConvex} \vline & \multicolumn{3}{c}{Ipopt} \vline & \multicolumn{3}{c}{Knitro}  \\
 & Inst & Obj & Gap  & Time  & Obj & Gap  & Time  & Obj & Gap  & Time  \\
\midrule
 & 10/01  & 0.31 & 2 &  -  & 0.32 &  -  & 1 & 0.32 &  -  & 0 \\
 & 10/02  & 0.22 & 3 &  -  & 0.22 &  -  & 1 & 0.22 &  -  & 0 \\
\parbox[t]{2mm}{\multirow{3}{*}{\rotatebox[origin=c]{90}{\texttt{ProductMode}}}} & 10/05  & 0.67 &  -  &  -  & 0.09 &  -  & 0 & 0.09 &  -  & 0 \\
 & 100/01  & 2.42 & 0 & 6 & 2.43 &  -  & 18 & 2.46 &  -  & 0 \\
 & 100/02  & 2.71 & 17 &  -  & 2.41 &  -  & 12 & 2.43 &  -  & 0 \\
 & 100/05  & 54.87 &  -  &  -  & 2.47 &  -  & 24 & 2.54 &  -  & 0 \\
 & 100/10  & 79.59 &  -  &  -  & 2.64 &  -  & 12 & 2.35 &  -  & 0 \\
 & 100/20  & 102.79 &  -  &  -  & 3.43 &  -  & 11 & 3.51 &  -  & 0 \\
 & 100/50  & 185.45 &  -  &  -  & 19.46 &  -  & 22 & 19.46 &  -  & 0 \\
 & 1000/01  & 25.21 & 0 &  -  & 25.13 &  -  & 272 & 25.10 &  -  & 1 \\
 & 1000/02  & 323.30 &  -  &  -  & 24.08 &  -  & 148 & 23.84 &  -  & 2 \\
 & 1000/05  & 533.37 &  -  &  -  & 25.70 &  -  & 121 & 24.89 &  -  & 4 \\
 \midrule
 & 10/01  & 0.30 & 1 &  -  & 0.32 &  -  & 0 & 0.33 &  -  & 0 \\
 & 10/02  & 0.25 & 17 &  -  & 0.22 &  -  & 0 & 0.22 &  -  & 0 \\
\parbox[t]{2mm}{\multirow{3}{*}{\rotatebox[origin=c]{90}{\texttt{StrongDualityMode}}}} & 10/05  & 0.13 &  -  &  -  & 0.09 &  -  & 0 & 0.09 &  -  & 0 \\
 & 100/01  & 2.42 & 0 & 0 & 2.44 &  -  & 2 & 2.43 &  -  & 0 \\
 & 100/02  & 2.41 & 4 &  -  & 2.44 &  -  & 2 & 2.53 &  -  & 0 \\
 & 100/05  & 2.33 & 7 &  -  & 2.50 &  -  & 0 & 2.32 &  -  & 0 \\
 & 100/10  & 2.41 & 39 &  -  & 2.21 &  -  & 1 & 2.08 &  -  & 0 \\
 & 100/20  & 3.55 & 145 &  -  & 3.43 &  -  & 8 & 2.90 &  -  & 0 \\
 & 100/50  & 64.20 &  -  &  -  & 23.34 &  -  & 1 & 185.45 &  -  & 0 \\
 & 1000/01  & 25.03 & 0 &  -  & 54.89 &  -  & 3 & 25.08 &  -  & 3 \\
 & 1000/02  & 23.80 & 0 &  -  & 71.32 &  -  & 9 & 23.78 &  -  & 6 \\
 & 1000/05  & 24.86 & 1 &  -  & 29.52 &  -  & 175 & 24.75 &  -  & 4 \\
\bottomrule
\end{tabular}
}
\caption{Gurobi NonConvex and NLP solvers with \texttt{ProductMode} and \texttt{StrongDualityMode}, Time in seconds (s), Gap in percent (\%).}
\label{table_prod}
\end{table}